\documentclass[12pt]{article}
\usepackage{latexsym,amssymb,amsmath,enumerate,float,geometry,cite,tikz}
\newtheorem{theorem}{Theorem}
\newtheorem{corollary}[theorem]{Corollary}

\newtheorem{lemma}[theorem]{Lemma}

\newtheorem{claim}{Claim}

\begin{document}

\tikzstyle{every node}=[circle, draw, inner sep=0pt, minimum width=4pt]

\title{\Large On some hard and some tractable cases of\\ the maximum acyclic matching problem}
\author{M. F\"{u}rst \and D. Rautenbach}
\date{}
\maketitle
\vspace{-10mm}
\begin{center}
{\small
Institute of Optimization and Operations Research, Ulm University,\\
Ulm, Germany, \texttt{maximilian.fuerst,dieter.rautenbach@uni-ulm.de}}
\end{center}

\begin{abstract}
Three well-studied types of subgraph-restricted matchings are 
induced matchings, 
uniquely restricted matchings, and
acyclic matchings.
While it is hard to determine the maximum size of a matching of each of these types,
whether some given graph has a maximum matching that is induced or 
has a maximum matching that is uniquely restricted, 
can both be decided efficiently.
In contrast to that we show that deciding whether a given bipartite graph of maximum degree at most four has a maximum matching that is acyclic is NP-complete.
Furthermore, we show that maximum weight acyclic matchings can be determined efficiently for $P_4$-free graphs and $2P_3$-free graphs,
and we characterize the graphs for which every maximum matching is acyclic.
\end{abstract}

{\small
\begin{tabular}{lp{12.5cm}}
\textbf{Keywords:} & 
Matching; 
induced matching; 
uniquely restricted matching;
acyclic matching
\end{tabular}
}

\pagebreak

\section{Introduction}

Three of the most natural types of subgraph-restricted matchings\cite{lopl,gohehela} are 
the induced matchings\cite{stva}, 
the uniquely restricted matchings\cite{gohile}, and 
the acyclic matchings \cite{gohehela}.
If $M$ is a matching in a graph $G$, and $G(M)$ is the subgraph of $G$ induced by the vertices of $G$ incident with some edge in $M$, then 
$M$ is {\it induced} if $G(M)$ is $1$-regular,
$M$ is {\it uniquely restricted} if $M$ is the only perfect matching of $G(M)$, and 
$M$ is {\it acyclic} if $G(M)$ is a forest.
Clearly, induced matchings are also acyclic.
Furthermore, since $M$ is uniquely restricted if and only if there is no $M$-alternating cycle \cite{gohile},
that is, a cycle whose every second edges belongs to $M$,
acyclic matchings are also uniquely restricted.
Therefore, 
if $\nu_s(G)$, $\nu_{ac}(G)$, and $\nu_{ur}(G)$ are the largest sizes of an induced, acyclic, and uniquely restricted matching in $G$, respectively,
then 
$$\nu_s(G)\leq \nu_{ac}(G)\leq \nu_{ur}(G)\leq \nu(G)$$
for every graph $G$,
where $\nu(G)$ is the classical matching number \cite{lopl},
that is, the largest size of an unrestricted matching.

Unlike the classical matching number, which is one of the most important tractable graph parameter, 
the three parameters $\nu_s(G)$, $\nu_{ac}(G)$, and $\nu_{ur}(G)$ are all computationally hard.
For the induced matching number $\nu_s(G)$, this has first been shown by Stockmeyer and Vazirani \cite{stva},
and strong inapproximability results are known \cite{dadelo,dumazi}.
For the uniquely restricted matching number $\nu_{ur}(G)$, the hardness was established by Golumbic, Hirst, and Lewenstein \cite{gohile},
and also for this parameter inapproximability results are known \cite{mi}.
The hardness of the acylic matching number $\nu_{ur}(G)$ was shown by Goddard, Hedetniemi, Hedetniemi, and Laskar \cite{gohehela}
by exploiting the well-known hardness of the maximum induced forest problem,
and Panda and Pradhan \cite{papr} provided further hardness results.

In view of these hardness results, it is quite surprising that 
\begin{itemize}
\item the graphs $G$ for which $\nu_s(G)=\nu(G)$ \cite{koro,cawa,dujoperaso,jora}, and also
\item the graphs $G$ for which $\nu_{ur}(G)=\nu(G)$ \cite{peraso}
\end{itemize}
can both be recognized efficiently.

As our first main result in this paper, we show that something similar does not hold for the acyclic matching number.
More precisely, we show that deciding $\nu_{ac}(G)=\nu(G)$ for a given graph $G$ is hard even under strong restrictions imposed on $G$.

Next to the hardness results mentioned above, 
some efficiently tractable cases have been studied 
for the induced matching number \cite{brmo,ca1,ca2,lo,koro,lomo},
the uniquely restricted matching number \cite{gohile,frja}, and 
the acyclic matching number \cite{bara, papr}.
We contribute to this line of research by showing that maximum weight acyclic matchings
can be found efficiently for $P_4$-free graphs and for $2P_3$-free graphs.
Finally, we characterize the graphs for which every maximum matching is acyclic.

The next section contains all our results and in a third section we discuss some open problems.

\medskip

\noindent Before we proceed to our results, we collect some terminology and notation.
We consider finite, simple, and undirected graphs.
Let $G$ be a graph.
The vertex set of $G$ is denoted $V(G)$ and the edge set of $G$ is denoted $E(G)$.
The neighborhood $N_G(u)$ of a vertex $u$ in $G$ is the set $\{ v\in V(G):uv\in E(G)\}$,
and the degree $d_G(u)$ of $u$ in $G$ is $|N_G(u)|$.
For a set $U$ of vertices of $G$, let $N_G(U)$ be $\bigcup_{u\in U}N_G(u)$.
A matching $M$ in a graph $G$ is a set of pairwise disjoint edges.
The set of vertices of $G$ incident with an edge in $M$ is $V(M)$,
and the subgraph of $G$ induced by $V(M)$ is $G(M)$.
For an integer $k$, let $[k]$ be the set of all positive integers at most $k$.

\section{Results}

We begin with our first main result showing that deciding the equality of the acyclic matching number and the matching number is hard 
even for quite restricted instances.

\begin{theorem}\label{theorem1}
For a given bipartite graph $G$ of maximum degree at most $4$, partite sets $A$ and $B$ with $|A|\leq |B|$, and matching number $|A|$, it is NP-complete to decide whether $\nu_{ac}(G)=\nu(G)$.
\end{theorem}
{\it Proof:} The proof relies on a restricted version of {\sc Satisfiability}.
It is known (cf.~[LO1] in \cite{gajo}) that {\sc Satisfiability} remains NP-complete when restricted to instances where every clause contains at most three literals, and for every variable $x$, at most three clauses contain the variable $x$ or its negation $\bar{x}$.
Note 
that variables that appear only positively or only negatively can be eliminated,
that clauses of size one allow to eliminate a variable,
and that clauses containing a variable and its negation can be removed.
Therefore, {\sc Satisfiability} remains NP-complete when restricted to instances 
where every clause contains two or three literals, and, 
for every variable $x$, 
at most two clauses contain the variable $x$, 
at most two clauses contain its negation $\bar{x}$, and
no clause contains both $x$ and $\bar{x}$.

Let $f$ be such an instance of {\sc Satisfiability} consisting of the clauses $c_1,\ldots,c_m$ over the 
boolean variables $x_1,\ldots,x_n$.
We construct a graph $G$ as in the statement of the theorem such that the order of $G$ is polynomially bounded in terms of $n$ and $m$, and $f$ is satisfiable if and only if $\nu_{ac}(G)=\nu(G)$. 

Starting with $G$ being the empty graph, 
\begin{itemize}
\item for every $i$ in $[n]$, we add to $G$ a subgraph $G_i$ as shown in Figure \ref{fig1}, and
\item for every $r$ in $[m]$, and every literal $\ell$ in the clause $c_r$, 
we add to $G$ a complete subgraph with the two vertices $\ell(r,1)$ and $\ell(r,2)$.
\end{itemize}
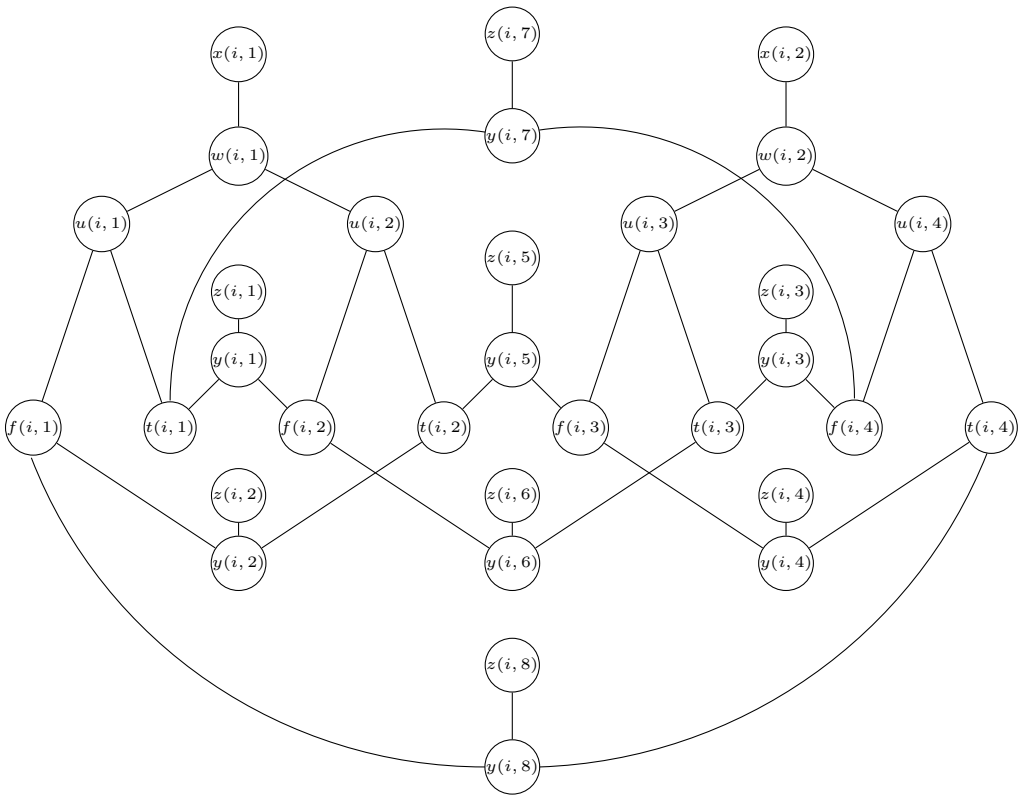
\begin{figure}[H]
\centering\tiny
\begin{tikzpicture}[scale = 0.9] 
	    \node (A) at (7,0) {$y(i,8)$}; 
	    \node (B) at (7,1.5) {$z(i,8)$};
	    \node (C) at (3,3) {$y(i,2)$}; 
	    \node (D) at (7,3) {$y(i,6)$}; 
	    \node (E) at (11,3) {$y(i,4)$}; 
	    \node (F) at (3,4) {$z(i,2)$};
	    \node (G) at (7,4) {$z(i,6)$}; 
	    \node (H) at (11,4) {$z(i,4)$};
	    \node (I) at (0,5) {$f(i,1)$};
	    \node (J) at (2,5) {$t(i,1)$}; 
	    \node (K) at (4,5) {$f(i,2)$}; 
	    \node (L) at (6,5) {$t(i,2)$};
	    \node (M) at (8,5) {$f(i,3)$}; 
	    \node (N) at (10,5) {$t(i,3)$}; 
	    \node (O) at (12,5) {$f(i,4)$};
	    \node (P) at (14,5) {$t(i,4)$}; 
	    \node (Q) at (3,6) {$y(i,1)$}; 
	    \node (R) at (7,6) {$y(i,5)$};
	    \node (S) at (11,6) {$y(i,3)$}; 
	    \node (T) at (3,7) {$z(i,1)$}; 
	    \node (U) at (7,7.5) {$z(i,5)$};
	    \node (V) at (11,7) {$z(i,3)$};
	    \node (W) at (1,8) {$u(i,1)$}; 
	    \node (X) at (5,8) {$u(i,2)$};
	    \node (Y) at (9,8) {$u(i,3)$}; 
	    \node (Z) at (13,8) {$u(i,4)$};
	    \node (a) at (3,9) {$w(i,1)$};
	    \node (b) at (7,9.3) {$y(i,7)$};
	    \node (c) at (11,9) {$w(i,2)$}; 
	    \node (d) at (3,10.5) {$x(i,1)$};
	    \node (e) at (7,10.8) {$z(i,7)$}; 
	    \node (f) at (11,10.5) {$x(i,2)$};
   
	    \draw[-] (A) to (B);
	    \draw (7.4,0) arc(-88:-21.5 :7.3); %A to P
	    \draw (6.6,0) arc(-90:-159:7.1); %A to I
	    \draw[-] (C) to (F);
	    \draw[-] (C) to (I);
	    \draw[-] (C) to (L);
	    \draw[-] (D) to (G);
	    \draw[-] (D) to (K);
	    \draw[-] (D) to (N);
	    \draw[-] (E) to (H);
	    \draw[-] (E) to (M);
	    \draw[-] (E) to (P);
	    \draw[-] (I) to (W);
	    \draw[-] (J) to (W);
	    \draw[-] (J) to (Q);
	    \draw (2,5.4) arc(180:81.3:4); % J to b
	    \draw[-] (K) to (Q);
	    \draw[-] (K) to (X);
	    \draw[-] (L) to (X);
	    \draw[-] (L) to (R);
	    \draw[-] (M) to (R);
	    \draw[-] (M) to (Y);
	    \draw[-] (N) to (Y);
	    \draw[-] (N) to (S);
	    \draw[-] (O) to (S);
	    \draw (12,5.43) arc(0:98.7:4); % O to b
	    \draw[-] (O) to (Z);
	    \draw[-] (P) to (Z);
	    \draw[-] (W) to (a);
	    \draw[-] (Q) to (T);
	    \draw[-] (X) to (a);
	    \draw[-] (R) to (U);
	    \draw[-] (S) to (V);
	    \draw[-] (Y) to (c);
	    \draw[-] (Z) to (c);
	    \draw[-] (a) to (d);
	    \draw[-] (b) to (e);
	    \draw[-] (c) to (f);
    \end{tikzpicture}
\caption{The gadget $G_i$ corresponding to the boolean variable $x_i$.} \label{fig1}
\end{figure}
\noindent Note that the clause $c_r$ of the form $x_i\vee \bar{x}_j$ leads to a subgraph with 
four vertices $x_i(r,1)$, $x_i(r,2)$, $\bar{x}_j(r,1)$, and $\bar{x}_j(r,2)$, and 
two edges $x_i(r,1)x_i(r,2)$ and $\bar{x}_j(r,1)\bar{x}_j(r,2)$.
Similarly, a clause containing three literals leads to a subgraph with six vertices and three edges.
Therefore, the graph constructed so far has order at most $32n+6m$.

In order to complete the construction of $G$,
for every $i$ in $[n]$, 
we add to $G$ some further edges between the subgraph $G_i$ and the subgraphs corresponding to the clauses containing $x_i$ or $\bar{x}_i$.
Therefore, let $i\in [n]$.
Let $c_r$ and $c_s$ be the clauses that contain $x_i$, where $r\leq s$,
that is, if $x_i$ is contained in just one clause, then let $r=s$.
First, suppose that $c_r$ contains only two literals $x_i$ and $\ell$.
Let $j\in [n]\setminus \{ i\}$ be such that $\ell$ is in $\{ x_j,\bar{x}_j\}$.
If $i<j$, then add to $G$ the two edges $x_i(r,1) f(i,1)$ and $\ell(r,1)f(i,2)$, and,
if $i>j$, then add to $G$ the two edges $x_i(r,1) f(i,2)$ and $\ell(r,1)f(i,1)$.
Next, suppose that $c_r$ contains three literals $x_i$, $\ell$, and $\ell'$.
Let $j,k\in [n]\setminus \{ i\}$ be such that $\ell$ is in $\{ x_j,\bar{x}_j\}$
and $\ell'$ is in $\{ x_k,\bar{x}_k\}$, where we may assume that $j<k$.
In this case, add to $G$ the two edges $\ell(r,1) f(i,1)$ and $\ell'(r,1)f(i,2)$.
If $r<s$, that is, $x_i$ is contained in two different clauses,
then proceed similarly for the clause $c_s$ 
using the vertices $f(i,3)$ and $f(i,4)$ instead of the vertices $f(i,1)$ and $f(i,2)$.
Furthermore, proceed similarly for the one or two clauses that contain $\bar{x}_i$
using the vertices $t(i,1)$, $t(i,2)$, $t(i,3)$, and $t(i,4)$
instead of the vertices $f(i,1)$, $f(i,2)$, $f(i,3)$, and $f(i,4)$.

This completes the construction of $G$, and we continue with a series of claims.

The next claim establishes some structural properties of $G$.

\begin{claim}\label{claim1}
$G$ is bipartite, has maximum degree $4$, 
partite sets $A$ and $B$ with $|A|\leq |B|$,
and matching number $|A|$.
\end{claim}
{\it Proof of Claim \ref{claim1}:}
It follows immediate from the construction that $G$ has maximum degree of $4$.

For the two sets
$$
A=
\bigcup\limits_{i\in [n]}
\Big(\{ u(i,j):j\in [4]\}\cup \{ x(i,j):j\in [2]\}\cup \{ y(i,j):j\in [8]\}\Big)
\cup\bigcup_{j\in[m]}\,\,\,\,\bigcup_{\ell\,\,is\,\,a\,\,literal\,\,in\,\,c_j}\ell(j,1)
$$
and $B=V(G)\setminus A$,
it is easy to verify that $G$ is bipartite with partite sets $A$ and $B$.
The matching $M$ that contains all edges incident with endvertices
as well as,
for every $i$ in $[n]$, 
the four edges $u(i,j)t(i,j)$ with $j\in [4]$,
satisfies $A\subseteq V(M)$,
which completes the proof of the claim. $\Box$

\medskip 

\noindent The next claim restricts the structure of maximum matchings in $G$
that are acyclic. 

\begin{claim}\label{claim2}
Every maximum matching $M$ in $G$ that is acyclic 
contains all edges incident with endvertices,
and, for every $i$ in $[n]$, the set $V(M)\cap V(G_i)$ is 
either
$V(G_i)\setminus \{ t(i,j):j\in [4]\}$
or
$V(G_i)\setminus \{ f(i,j):j\in [4]\}$.
\end{claim}
{\it Proof of Claim \ref{claim2}:} 
By Claim \ref{claim1}, we have $A\subseteq V(M)$.
Let $i$ be in $[n]$.
Since $x(i,j)\in A$ for $j\in [2]$,
the matching $M$ contains the two edges incident with the endvertices $x(i,1)$ and $x(i,2)$.

If $u(i,1)t(i,1),u(i,2)f(i,2)\in M$,
then $G(M)$ contains the cycle 
$$u(i,1)t(i,1)y(i,1)f(i,2)u(i,2)w(i,1),$$
which is a contradiction.
Hence, by symmetry, the matching $M$ contains 
\begin{itemize}
\item either $u(i,j)t(i,j)\in M$ for $j\in [2]$ or $u(i,j)f(i,j)\in M$ for $j\in [2]$, and
\item either $u(i,j)t(i,j)\in M$ for $j\in [4]\setminus [2]$ or $u(i,j)f(i,j)\in M$ for $j\in [4]\setminus [2]$.
\end{itemize}
If $u(i,j)t(i,j)\in M$ for $j\in [2]$ and $u(i,j)f(i,j)\in M$ for $j\in [4]\setminus [2]$,
then $G(M)$ contains the cycle 
$$u(i,2)t(i,2)y(i,5)f(i,3)u(i,3)w(i,2)u(i,4)f(i,4)y(i,7)t(i,1)u(i,1)w(i,1)u(i,2),$$
which is a contradiction. 
Hence, either $u(i,j)t(i,j)\in M$ for $j\in [4]$ or $u(i,j)f(i,j)\in M$ for $j\in [4]$.

If $f(i,1)\in V(M)$ and $u(i,j)t(i,j)\in M$ for $j\in [4]$, 
then $G(M)$ contains the cycle
$$f(i,1)y(i,2)t(i,2)u(i,2)w(i,1)u(i,1)f(i,1),$$
which is a contradiction.
Hence, by symmetry, it follows that the set $V(M)\cap V(G_i)$ is
either
$V(G_i)\setminus \{ t(i,j):j\in [4]\}$
or
$V(G_i)\setminus \{ f(i,j):j\in [4]\}$.
This implies that $y(i,j)z(i,j)\in M$ for $j\in [8]$,
and that $\ell(j,1)\ell(j,2)\in M$ for every $j\in [m]$ and every literal $\ell$ in the clause $c_j$,
which completes the proof of the claim.
$\Box$

\medskip

\noindent By construction, for every $i$ in $[n]$, the only vertices within $V(G_i)$ 
that can have neighbors outside of $V(G_i)$
are the eight vertices $t(i,j)$ and $f(i,j)$ for $j\in [4]$.
By Claim \ref{claim2}, if $M$ is a maximum matching in $G$ that is acyclic, then 
$G[V(M)\cap V(G_i)]$ is a forest.
Furthermore,
if $V(M)\cap V(G_i)=V(G_i)\setminus \{ t(i,j):j\in [4]\}$,
then the vertices $f(i,1)$ and $f(i,2)$ lie in one component of $G[V(M)\cap V(G_i)]$,
and the vertices $f(i,3)$ and $f(i,4)$ lie in another component of $G[V(M)\cap V(G_i)]$,
and,
if $V(M)\cap V(G_i)=V(G_i)\setminus \{ f(i,j):j\in [4]\}$,
then the vertices $t(i,1)$ and $t(i,2)$ lie in one component of $G[V(M)\cap V(G_i)]$,
and the vertices $t(i,3)$ and $t(i,4)$ lie in another component of $G[V(M)\cap V(G_i)]$. 

\begin{claim}\label{claim3}
$\nu_{ac}(G)=\nu(G)$ if and only if $f$ is satisfiable.
\end{claim}
{\it Proof of Claim \ref{claim3}:} 
First, suppose that $\nu_{ac}(G)=\nu(G)$.
Let $M$ be a maximum matching in $G$ that is acyclic.
For every $i$ in $[n]$, set the variable $x_i$ to \texttt{true}
if and only if $t(i,1)\in V(M)$.
Now, let $r$ be in $[m]$.
First, suppose that $c_r$ contains two literals 
$\ell\in \{ x_i,\bar{x}_i\}$ and
$\ell'\in \{ x_j,\bar{x}_j\}$ for some $i$ and $j$ in $[n]$.
If both literals in $c_r$ are not \texttt{true},
then $G(M)$ contains 
a path between $\ell(r,1)$ and $\ell'(r,1)$ through $G_i$
and also through $G_j$,
which contradicts the assumption that $M$ is acyclic.
Next, suppose that $c_r$ contains three literals
$\ell\in \{ x_i,\bar{x}_i\}$,
$\ell'\in \{ x_j,\bar{x}_j\}$, and
$\ell''\in \{ x_k,\bar{x}_k\}$ for some $i$, $j$, and $k$ in $[n]$.
If all three literals in $c_r$ are not \texttt{true},
then $G(M)$ contains 
a path between $\ell(r,1)$ and $\ell'(r,1)$ through $G_k$,
a path between $\ell'(r,1)$ and $\ell''(r,1)$ through $G_i$, and
a path between $\ell(r,1)$ and $\ell''(r,1)$ through $G_j$,
which contradicts the assumption that $M$ is acyclic.
Hence, the considered truth assignment satisfies every clause of $f$.

Next, suppose that $f$ is satisfiable, and consider a satisfying truth assignment.
Let $M$ be the maximum matching in $G$ that contains all edges incident with endvertices,
and satisfies
$V(M)\cap V(G_i)=V(G_i)\setminus \{ f(i,j):j\in [4]\}$
if and only if $x_i$ is \texttt{true} for every $i$ in $[n]$.
The existence of such a matching follows easily from the construction.
Note that, if $x_i$ is not \texttt{true},
then $V(M)\cap V(G_i)=V(G_i)\setminus \{ t(i,j):j\in [4]\}$.

Suppose for a contradiction, that $G(M)$ contains a cycle $C$.
It is easy to see that $C$ contains a vertex $\ell(r,1)$
for some $r$ in $[m]$ and some literal $\ell$ in the clause $c_r$.
If the clause $c_r$ contains only one further literal $\ell'$,
and $i$ and $j$ in $[n]$ are such that 
$\ell\in \{ x_i,\bar{x}_i\}$ and
$\ell'\in \{ x_j,\bar{x}_j\}$, 
then the two neighbors of $\ell(r,1)$ in $G_i$ and $G_j$ are both in $V(M)$.
By the construction of $G$ and the definition of $M$, 
this implies that both literals $\ell$ and $\ell'$ are not \texttt{true},
which is a contradiction.
Hence, we may assume that the clause $c_r$ 
contains two further literals $\ell'$ and $\ell''$.
Let $i$, $j$, and $k$ in $[n]$ be such that 
$\ell\in \{ x_i,\bar{x}_i\}$,
$\ell'\in \{ x_j,\bar{x}_j\}$, and
$\ell''\in \{ x_k,\bar{x}_k\}$.
Since $\ell(r,1)$ lies on $C$, 
the two neighbors of $\ell(r,1)$ in $G_j$ and $G_k$ are both in $V(M)$.
In view of the intersection of $V(M)$ with $V(G_j)$ and $V(G_k)$,
it follows easily that $C$ contains the vertices $\ell'(r,1)$ and $\ell''(r,1)$.
This implies that 
the two neighbors of $\ell'(r,1)$ in $G_i$ and $G_k$ are both in $V(M)$,
and that also
the two neighbors of $\ell''(r,1)$ in $G_i$ and $G_j$ are both in $V(M)$.
By the definition of $M$, 
this implies that all three literals $\ell$, $\ell'$, and $\ell''$ are not \texttt{true},
which is a contradiction, and completes the proof of the claim. 
$\Box$

\medskip

\noindent The above three claims clearly imply the desired statement,
which completes the proof. 
$\Box$

\medskip

\noindent Subdividing each non-pendant edge in the graph $G$ 
constucted within the proof of Theorem \ref{theorem1} 
an even number of times,
and attaching endvertices to the vertices created by these subdivisions
easily implies that the conclusion of Theorem \ref{theorem1} still holds
when we additionally require a large girth.
Furthermore, since $\nu(G)=\nu_{ur}(G)$ for the graph $G$
constucted within the proof of Theorem \ref{theorem1},
the proof of this result also implies the hardness of deciding 
whether $\nu_{ur}(G)$ equals $\nu_{ac}(G)$ 
for a given graph $G$ as in Theorem \ref{theorem1}.

The following three results concern tractable cases for the maximum weight acyclic matching problem.

\begin{theorem}\label{theorem2}
For a given edge-weighted $P_4$-free graph, 
a maximum weight acyclic matching can be determined in polynomial time.
\end{theorem}
{\it Proof:} Is it well-known that every $P_4$-free graph $G$ of order at least $2$
is either disconnected 
or the join of two disjoint $P_4$-free graphs.
Clearly, a maximum weight acyclic matching of a disconnected graph is the union of maximum weight acyclic matchings of its components.
Furthermore, if $G$ is the join of $G_1$ and $G_2$,
and $M$ is an acyclic matching in $G$, then 
either $M\subseteq E(G_1)$,
or $M\subseteq E(G_2)$,
or $M$ consists of exactly one edge between a vertex of $G_1$ and a vertex of $G_2$.
These observations lead to a simple polynomial time reduction algorithm 
that determines a maximum weight acyclic matching 
in a given edge-weighted $P_4$-free graph.
$\Box$

\medskip

\noindent Our next lemma shows that for $2P_3$-free graphs, 
induced matchings and acyclic matching are similar.

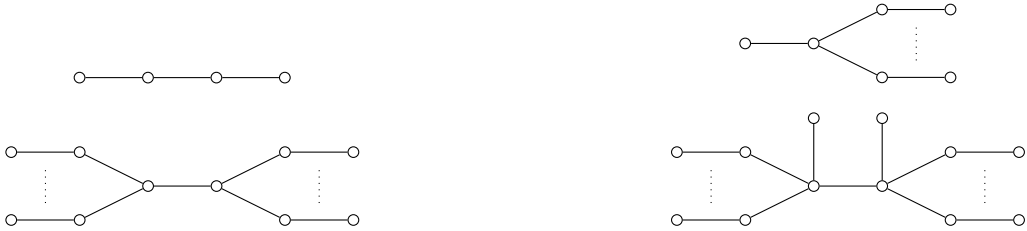
\begin{figure}[H]
\begin{minipage}[t]{0.48\textwidth}
\centering\tiny
\begin{tikzpicture}[scale = 0.9] 
	    \node (A) at (0,0) {};
	    \node (B) at (1,0) {};
	    \node (C) at (2,0) {};
	    \node (D) at (3,0) {};
	    
	    \draw[-] (A) to (B);
	    \draw[-] (B) to (C);
	    \draw[-] (C) to (D);
\end{tikzpicture}
\end{minipage} \hfill
\begin{minipage}[t]{0.49\textwidth}
\centering\tiny
\begin{tikzpicture}[scale = 0.9] 
	    \node (A) at (0,0) {};
	    \node (B) at (1,0) {};
	    \node (C) at (2,0.5) {};
	    \node (D) at (3,0.5) {};
	    \node (E) at (2,-0.5) {};
	    \node (F) at (3,-0.5) {};
	    
	    \draw[-] (A) to (B) to (C) to (D);
	    \draw[-] (B) to (E) to (F);
	    \draw[dotted] (2.5,-0.25) to (2.5,0.25);
\end{tikzpicture}
\end{minipage} \\[10pt]
\begin{minipage}[t]{0.48\textwidth}
\centering\tiny
\begin{tikzpicture}[scale = 0.9] 
	    \node (A) at (0,0.5) {};
	    \node (B) at (1,0.5) {};
	    \node (C) at (2,1) {};
	    \node (D) at (3,1) {};
	    \node (E) at (2,0) {};
	    \node (F) at (3,0) {};
	    \node (G) at (-1,1) {};
	    \node (H) at (-2,1) {};
	    \node (I) at (-1,0) {};
	    \node (J) at (-2,0) {};
	    
	    \draw[-] (A) to (B) to (C) to (D);
	    \draw[-] (B) to (E) to (F);
	    \draw[-] (A) to (G) to (H);
	    \draw[-] (A) to (I) to (J);
	    \draw[dotted] (2.5,0.25) to (2.5,0.75);
	    \draw[dotted] (-1.5,0.25) to (-1.5,0.75);
\end{tikzpicture}
\end{minipage} \hfill
\begin{minipage}[t]{0.49\textwidth}
\centering\tiny
 \begin{tikzpicture}[scale = 0.9] 
	    \node (A) at (0,0.5) {};
	    \node (B) at (1,0.5) {};
	    \node (C) at (2,1) {};
	    \node (D) at (3,1) {};
	    \node (E) at (2,0) {};
	    \node (F) at (3,0) {};
	    \node (G) at (-1,1) {};
	    \node (H) at (-2,1) {};
	    \node (I) at (-1,0) {};
	    \node (J) at (-2,0) {};
	    \node (K) at (0,1.5) {};
	    \node (L) at (1,1.5) {};
	    
	    \draw[-] (A) to (B) to (C) to (D);
	    \draw[-] (B) to (E) to (F);
	    \draw[-] (A) to (G) to (H);
	    \draw[-] (A) to (I) to (J);
	    \draw[-] (A) to (K);
	    \draw[-] (B) to (L);
	    \draw[dotted] (2.5,0.25) to (2.5,0.75);
	    \draw[dotted] (-1.5,0.25) to (-1.5,0.75);
 \end{tikzpicture}
\end{minipage}
\caption{An illustration of Lemma \ref{lemma1} $(i)-(iv)$.} \label{fig2}
\end{figure}

\begin{lemma}\label{lemma1}
If $M$ is an acyclic matching in a $2P_3$-free graph $G$, 
then $G(M)$ has at most one component that is not a $K_2$,
and such a component $T$ satisfies 
\begin{enumerate}[(i)]
\item either $T$ is a $P_4$,
\item or 
\begin{eqnarray*}
V(T) &=& \{ x,x'\}\cup \{ y_i:i\in [k]\}\cup \{ y'_i:i\in [k]\}\mbox{ and}\\
E(T)&=& \{ xx'\}\cup \{ y_iy_i':i\in [k]\}\cup \{ xy_i:i\in [k]\}
\end{eqnarray*}
for some integer $k$ at least $2$,
\item or 
\begin{eqnarray*}
V(T) &=& \{ x,x'\}\cup \{ y_i:i\in [k]\}\cup \{ y'_i:i\in [k]\}\cup \{ z_i:i\in [\ell]\}\cup \{ z'_i:i\in [\ell]\}\mbox{ and}\\
E(T)&=& \{ xx'\}\cup \{ y_iy_i':i\in [k]\}\cup \{ xy_i:i\in [k]\}\cup \{ z_iz_i':i\in [\ell]\}\cup \{ x'z_i:i\in [\ell]\}
\end{eqnarray*}
for some positive integers $k$ and $\ell$,
\item or 
\begin{eqnarray*}
V(T) &=& \{ x,x',y,y'\}\cup \{ w_i:i\in [k]\}\cup \{ w'_i:i\in [k]\}\cup \{ z_i:i\in [\ell]\}\cup \{ z'_i:i\in [\ell]\}\mbox{ and}\\
E(T)&=& \{ xx',yy',xy\}\cup \{ w_iw_i':i\in [k]\}\cup \{ xw_i:i\in [k]\}\cup \{ z_iz_i':i\in [\ell]\}\cup \{ yz_i:i\in [\ell]\}
\end{eqnarray*}
for some positive integers $k$ and $\ell$.
\end{enumerate}
\end{lemma}
{\it Proof:} Let $F=G(M)$. 
Since $M$ is acyclic, $F$ is a forest with a perfect matching.
Since every component of $F$ that is not a $K_2$, contains a $P_3$,
the $2P_3$-freeness of $G$ implies that at most one component of $F$ is not a $K_2$.
Let $T$ be such a component.
Let $P:u_1\ldots u_n$ be a longest path in $T$.
Since $F$ has a perfect matching, we obtain $n\geq 4$, and,
since $G$ is $2P_3$-free, we obtain $n\leq 6$.
Note that $d_F(u_1)=d_F(u_n)=1$, 
which implies that $M$ contains the two edges $u_1u_2$ and $u_{n-1}u_n$.
By the choice of $P$, we obtain furthermore that $N_F(u_2)=\{ u_1,u_3\}$ and $N_F(u_{n-1})=\{ u_{n-2},u_n\}$.
If $n=4$, this already implies that $T$ is a $P_4$.
For $n=5$, the matching $M$ contains an edge $u_3u_3'$ for some $u_3\in V(T)\setminus V(P)$.
By symmetry and the choice of $P$, 
it follows that all vertices in $N_{T-(E(T)\cap M)}(u_3)=N_T(u_3)\setminus \{ u_3'\}$ have degree $2$ in $T$,
and that their neighbors distinct from $u_3$ are endvertices of $T$.
Altogether, it follows that $T$ is as in (ii).
Now, let $n=6$.
Since $M$ contains a perfect matching of $T$,
it follows that either $u_3u_4\in M$
or $M$ contains two edges $u_3u_3'$ and $u_4u_4'$ for some distinct vertices $u_3'$ and $u_4'$ in $V(T)\setminus V(P)$.
As before it follows that all vertices in 
$\left(N_{T-(E(T)\cap M)}(u_3)\cup N_{T-(E(T)\cap M)}(u_4)\right)\setminus \{ u_3,u_4\}$ have degree $2$ in $T$,
and that their neighbors distinct from $u_3$ and $u_4$ are endvertices of $T$.
Hence, if $u_3u_4\in M$, then $T$ is as in (iii), and, otherwise, $T$ is as in (iv).
$\Box$

\begin{theorem}\label{theorem3}
For a given edge-weighted $2P_3$-free graph, 
a maximum weight acyclic matching can be determined in polynomial time.
\end{theorem}
{\it Proof:} Let $G$ be an edge-weighted $2P_3$-free graph of order $n$. 
Using Lemma \ref{lemma1}, we will reduce the problem to determine a maximum weight acyclic matching in $G$
to polynomially many instances of the maximum weight induced matching problem in smaller $2P_3$-free graphs.
As shown by Lozin and Mosca \cite{lomo}
a maximum weight induced matching can be determined in polynomial time for a given $2P_3$-free graph.

By Lemma \ref{lemma1}, a maximum weight acyclic matching $M$ in $G$ is a matching of largest weight of one of the following five types:
\begin{enumerate}[(1)]
\item An induced matching.
\item An acyclic matching $M$ for which $G(M)$ contains a component $T$ that is a $P_4$.
\item An acyclic matching $M$ for which $G(M)$ contains a component $T$ as in (ii) of Lemma \ref{lemma1}.
\item An acyclic matching $M$ for which $G(M)$ contains a component $T$ as in (iii) of Lemma \ref{lemma1}.
\item An acyclic matching $M$ for which $G(M)$ contains a component $T$ as in (iv) of Lemma \ref{lemma1}.
\end{enumerate}
We already observed that a maximum weight matching of type (1) 
can be determined efficiently \cite{lomo}.
A maximum weight matching $M$ of type (2) 
can be determined by considering all $O(n^4)$ choices for the $P_4$ component $T$,
adding to $M$ the two edges of a perfect matching of $T$
as well as a maximum weight induced matching of the graph $G-N_G[V(T)]$.

A maximum weight matching $M$ of type (3) 
can be determined by 
\begin{itemize}
\item considering all $O(n^6)$ choices for the six vertices 
$x$, $x'$, $y_1$, $y_1'$, $y_2$, and $y_2'$ of $T$,
\item adding to $M$ the three edges $xx'$, $y_1y_1'$, and $y_2y_2'$,
as well as 
\item a maximum weight induced matching $M'$ of the graph $G'=G-N_G[\{ x',y_1,y_1',y_2,y_2'\}]$,
where the weight of the edges of $G'$ that are between vertices in $N_G(x)$ is changed to $-1$.
\end{itemize}
By the construction of $G'$ and the change of the weight function, every edge in $M'$ contains at most one neighbor of $x$.
Note that the edges in $M'$ containing a neighbor of $x$, 
correspond to the edges
$y_3y_3',\ldots,y_ky_k'$ as in Lemma \ref{lemma1}(ii).

In the remaining cases, we proceed similarly.

A maximum weight matching $M$ of type (4) 
can be determined by 
considering all $O(n^6)$ choices for the six vertices $x$, $x'$, $y_1$, $y_1'$, $z_1$, and $z_1'$ of $T$,
adding to $M$ the three edges $xx'$, $y_1y_1'$, and $z_1z_1'$,
as well as a maximum weight induced matching $M'$ of the graph $G'=G-\left(N_G[\{ y_1,y_1',z_1,z_1'\}]\cup \left(N_G(x)\cap N_G(x')\right)\right)$,
where the weight of the edges of $G'$ that are between vertices in $N_G(x)\cup N_G(x')$ is changed to $-1$.
Note that 
$N_G(u)\cap \{ x,x',y_1,y_1',z_1,z_1'\}$
is either $\emptyset$, or $\{ x\}$, or $\{ x'\}$
for every vertex $u$ in $V(G')$,
and that every edge in $M'$ contains at most one vertex with a neighbor in $\{ x,x'\}$,
which implies that the structure of $G(M)$ is as desired.

Finally,
a maximum weight matching $M$ of type (5) can be determined 
by considering all $O(n^8)$ choices for the eight vertices $x$, $x'$, $y$, $y'$, $w_1$, $w_1'$, $z_1$, and $z_1'$ of $T$,
adding to $M$ the four edges $xx'$, $yy'$, $w_1w_1'$, and $z_1z_1'$,
as well as a maximum weight induced matching of the graph $G'=G-\left(N_G[\{ x',y',w_1,w_1',z_1,z_1'\}]\cup \left(N_G(x)\cap N_G(y)\right)\right)$,
where the weight of the edges of $G'$ that are between vertices in $N_G(x)\cup N_G(y)$ is changed to $-1$.

Altogether, we need to solve $O(n^8)$ instances of the maximum weight induced matching problem in edge-weighted $2P_3$-free graphs that are induced subgraphs of $G$ such that the encoding lengths of the weight functions are polynomially bounded in terms of the encoding length of the original weight function, which, by \cite{lomo}, completes the proof. $\Box$

\medskip

\noindent Our final result concerns the graphs for which every maximum matching is acyclic.

\begin{theorem}\label{theorem4}
Let $G$ be a graph.

Every maximum matching in $G$ is acyclic
if and only if 
every component of $G$ is a tree or an odd cycle.
\end{theorem}
{\it Proof:} Clearly, if every component of $G$ is a tree or an odd cycle,
then every maximum matching in $G$ is acyclic.
Now, let $G$ be such that every maximum matching in $G$ is acyclic.
Let $C$ be a cycle in $G$.
In order to complete the proof, 
it suffices to show that $C$ is a component of $G$ of odd order.
Let $M$ be a maximum matching in $G$ maximizing $|V(M)\cap V(C)|$.
Since $M$ is acyclic, not all vertices of $C$ are covered by $M$.

Suppose that $M$ contains an edge between a vertex on $C$ and a vertex not on $C$.
If $u_0\ldots u_k$ is a shortest path within $C$ between a vertex $u_0$ not covered by $M$
and a vertex $u_k$ such that $u_kv\in M$ for some vertex $v$ not on $C$,
then $k$ is odd, and $u_{2i-1}u_{2i}\in M$ for every $i\in [(k-1)/2]$.
Now,
$$M'=\left(M\setminus \left(\{ u_kv\}\cup \left\{ u_{2i-1}u_{2i}:i\in [(k-1)/2]\right\}\right)\right)\cup \left\{ u_{2i-2}u_{2i-1}:i\in [(k+1)/2]\right\}$$
is a maximum matching in $G$ with $|V(M')\cap V(C)|>|V(M)\cap V(C)|$,
contradicting the choice of $M$.
Hence, $M$ contains no edge between a vertex on $C$ and a vertex not on $C$.
Since $M$ is a maximum matching, 
this easily implies that $M$ covers all but exactly one vertex of $C$,
which implies that $C$ has odd order.

Suppose that $C$ is not a component of $G$.
This implies that some vertex $u$ on $C$ has a neighbor $v$ that is not on $C$.
Removing from $M$ all edges that belong to $C$ as well as any edge incident with $v$,
and adding to the resulting matching 
the edge $uv$
as well as a matching in $C$ covering all vertices of $C$ except for $u$,
results in a maximum matching in $G$ that is not acyclic,
which is a contradiction.
Hence, $C$ is a component of $G$, which completes the proof. $\Box$

\begin{corollary}\label{corollary1}
Let $G$ be a graph.

Every maximum matching in $G$ is induced
if and only if 
every component of $G$ is a star or a triangle.
\end{corollary}
{\it Proof:} Clearly, if every component of $G$ is a star or a triangle,
then every maximum matching is induced.
Now, let $G$ be such that every maximum matching in $G$ is induced.
Since this implies that every maximum matching in $G$ is acyclic,
the graph $G$ is as in Theorem \ref{theorem4}.
Since odd cycles of length at least $5$ 
have non-induced maximum matchings,
every component of $G$ that is not a tree, is a triangle.
Now, suppose that $T$ is a component of $G$ that is a tree but not a star.
Let $M$ be a maximum matching of $G$ that contains a pendant edge $uv$ of $T$,
where $u$ is an endvertex.
Since $T$ is not a star, the vertex $v$ has a neighbor $w$ that is not an endvertex.
Since $M$ is a maximum matching that is induced, the vertex $w$ is not covered by $M$
but every vertex in $N_T(w)$ is covered by $M$.
Removing an edge from $M$ incident with a neighbor $x$ of $w$ that is distinct from $v$, and adding $wx$ to the resulting matching,
yields a non-induced maximum matching, which is a contradiction. $\Box$

\section{Conclusion}

Two obvious open problems motivated by our results concern 
\begin{itemize}
\item the complexity of deciding whether $\nu_{ac}(G)=\nu(G)$ 
for a given bipartite graph $G$ of maximum degree $3$, and 
\item the complexity of the maximum acyclic matching problem for $P_5$-free graphs.
\end{itemize}
We believe that the first of these two problems can be solved efficiently,
and add some comments concerning the second problem.
Kobler and Rotics \cite{koro} describe efficient algorithms for the maximum induced matching problem restricted to some subclasses of $P_5$-free graphs. 
Their key insight is that the square of the line graph $L^2(G)$ of a $P_5$-free graph $G$ is still $P_5$-free, and that induced matchings in $G$ 
correspond to independent sets in $L^2(G)$.
In view of \cite{lovavi}, it actually follows that the maximum induced matching problem can be solved efficiently for all $P_5$-free graphs.
If $M$ is an acyclic matching in a $P_5$-free graph $G$, then $G(M)$ is the disjoint union of $P_4$s and $K_2$s, and, may, in particular, contain many $P_4$s. Therefore, for the maximum acyclic matching problem in $P_5$-free graphs, an approach as in the proof of Theorem \ref{theorem3} does not seem to work. Nevertheless, considering $M$ as a vertex subset of $L^2(G)$, there might be a connection with the dissociation set problem.

Note that if $M$ is a maximum matching in a graph $G$ that is acyclic,
then $V(G)\setminus V(M)$ is a so-called independent feedback vertex set.
The algorithmic problems in connection with (minimum) independent feedback vertex sets 
have attracted considerable attention recently \cite{aggusash,bodafejopa,miphrasa,taitzh}.
Maybe some of the techniques from this context
lead to interesting cases, 
where deciding the equality of the acyclic matching number and the matching number is tractable.


\begin{thebibliography}{}
\bibitem{aggusash} A. Agrawal, S. Gupta, S. Saurabh, R. Sharma, Improved Algorithms and Combinatorial Bounds for Independent Feedback Vertex Set, Leibniz International Proceedings in Informatics 63 (2017) 2:1-2:14.
\bibitem{bara} J. Baste, D. Rautenbach, Degenerate Matchings and Edge Colorings, arXiv:1702.02358.
\bibitem{bodafejopa} M. Bonamy, K.K. Dabrowski, C. Feghali, M. Johnson, D. Paulusma, Independent Feedback Vertex Set for $P_5$-free Graphs, arXiv:1707.09402.
\bibitem{brmo} A. Brandst\"{a}dt, R. Mosca, On distance-$3$ matchings and induced matchings, Discrete Applied Mathematics 159 (2011) 509-520.
\bibitem{ca1} K. Cameron, Induced matchings, Discrete Applied Mathematics 24 (1989) 97-102.
\bibitem{ca2} K. Cameron, Induced matchings in intersection graphs, Discrete Mathematics 278 (2004) 1-9.
\bibitem{cawa} K. Cameron, T. Walker, The graphs with maximum induced matching and maximum matching the same size, Discrete Mathematics 299 (2005) 49-55.
\bibitem{dadelo} K.K. Dabrowski, M. Demange, V.V. Lozin, New results on maximum induced matchings in bipartite graphs and beyond,  Theoretical Computer Science 478 (2013) 33-40.
\bibitem{dujoperaso} M. Duarte, F. Joos, L.D. Penso, D. Rautenbach, U. Souza, Maximum induced matchings close to maximum matchings, Theoretical Computer Science 588 (2015) 131-137.
\bibitem{dumazi} W. Duckworth, D.F. Manlove, M. Zito, On the approximability of the maximum induced matching problem, Journal of Discrete Algorithms 3 (2005) 79-91.
\bibitem{frja} M.C. Francis, D. Jacob, S. Jana, Uniquely restricted matchings in interval graphs, arXiv:1604.07016.
\bibitem{gajo} M.R. Garey, D.S. Johnson, Computers and Intractability: A Guide to the Theory of NP-Completeness, W.H. Freeman and Co. pp. x+338, 1979. 
\bibitem{gohehela} W. Goddard, S.M. Hedetniemi, S.T. Hedetniemi, R. Laskar, Generalized subgraph-restricted matchings in graphs, Discrete Mathematics 293 (2005) 129-138.
\bibitem{gohile} M.C. Golumbic, T. Hirst, M. Lewenstein, Uniquely restricted matchings, Algorithmica 31 (2001) 139-154.
\bibitem{jora} F. Joos, D. Rautenbach, Equality of Distance Packing Numbers, Discrete Mathematics 338 (2015) 2374-2377.
\bibitem{koro} D. Kobler, U. Rotics, Finding Maximum Induced Matchings in Subclasses of Claw-Free and $P_5$-Free Graphs, and in Graphs with Matching and Induced Matching of Equal Maximum Size, Algorithmica 37 (2003) 327-346.
\bibitem{lovavi} D. Lokshantov, M. Vatshelle, Y. Villanger, Independent Set in $P_5$-free Graphs in Polynomial Time, Proceedings of the Twenty-fifth Annual ACM-SIAM Symposium on Discrete Algorithms (SODA '14), 570-581.
\bibitem{lo} V. Lozin, On maximum induced matchings in bipartite graphs, Information Processessing Letters 81 (2002) 7-11.
\bibitem{lomo} V. Lozin, R. Mosca, Maximum regular induced subgraphs in $2P_3$-free graphs, Theoretical Computer Science 460 (2012) 26-33.
\bibitem{lopl} L. Lov\'{a}sz, M.D. Plummer, Matching theory, Annals of Discrete Mathematics, 29, 1986, North-Holland Publishing Co., Amsterdam.
\bibitem{mi} S. Mishra, On the maximum uniquely restricted matching for bipartite graphs, Electronic Notes in Discrete Mathematics 37 (2011) 345-350.
\bibitem{miphrasa} N. Misra, G. Philip, V. Raman, S. Saurabh, On parameterized independent feedback vertex set, Theoretical Computer Science 461 (2012) 65-75.
\bibitem{papr} B.S. Panda, D. Pradhan, Acyclic matchings in subclasses of bipartite graphs, Discrete Mathematics, Algorithms and Applications 4 (2012) 1250050 (15 pages).
\bibitem{peraso} L.D. Penso, D. Rautenbach, U. Souza, Graphs in which some and every maximum matching is uniquely restricted, arXiv:1504.02250.
\bibitem{stva} L.J. Stockmeyer, V.V. Vazirani, NP-completeness of some generalizations of the maximum matching problem, Information Processing Letters 15 (1982) 14-19.
\bibitem{taitzh} Y. Tamura, T. Ito, X. Zhou, Algorithms for the independent feedback vertex set problem, IEICE Transactions on Fundamentals of Electronics, Communications and Computer Sciences 98 (2015) 1179-1188.
\end{thebibliography}
\end{document}